# On Chromatic no. of $3K_1$-free graphs and R(3, k)

Medha Dhurandhar

**Abstract**

Here we prove that if G is $\{3K_1, K_{\omega+1}\}$-free with $\omega \leq 11$, then (1) $\chi \leq \frac{1}{8}(\omega^2+12\omega-13)$, if $\omega$ is odd, and (2) $\chi \leq \frac{1}{8}(\omega^2+10\omega)$, if $\omega$ is even. We further conjecture that the results are true in general for all $\omega$. We also conjecture that (A) if $\omega$ is odd and $R(3, \omega)$ is even, then $R(3, \omega) = \frac{1}{4}(\omega^2+8\omega-9)$, (B) if $\omega$ and $R(3, \omega)$ are both odd, then $\frac{1}{4}(\omega^2+8\omega-13)$, (C) if $\omega$ and $R(3, \omega)$ are both even, then $R(3, \omega) = \frac{1}{4}(\omega^2+6\omega)$ and (D) if $\omega$ is even and $R(3, \omega)$ is odd, then $\frac{1}{4}(\omega^2+6\omega-4)$. Again we verify the results for $\omega \leq 9$.

## 1. Introduction

We consider in this paper only finite, simple, connected, undirected graphs. The vertex set of G is denoted by V(G), the edge set by E(G), the set of neighbors of v in G by N(v), the maximum independence number by $\alpha(G)$, the maximum degree of vertices in G by $\Delta(G)$, the maximum clique size by $\omega(G)$ and the chromatic number by $\chi(G)$.

In 1998, Reed [12] conjectured that $\chi(G) \leq \lceil \frac{1}{2}(\Delta+\omega+1) \rceil$. Please refer to [1, 2, 3, 4, 5, 6, 7, 8, 9] for some classes of graphs for which Reed's conjecture holds good.. In [2], it was proved that $\chi(G) \leq \lceil(\Delta+\omega+1)\rceil$ if $\alpha(G) = 2$. Also using data from [12], it was proved in [13] that if G is $\{3K_1, K_5\}$-free, then $\chi(G) \leq 7$ and if G is $\{3K_1, K_6\}$-free, then $\chi(G) \leq 9$. Here we prove that if G is $\{3K_1, K_{\omega+1}\}$-free with $\omega \leq 11$, then (1) $\chi \leq \frac{1}{8}(\omega^2+12\omega-13)$, if $\omega$ is odd, and (2) $\chi \leq \frac{1}{8}(\omega^2+10\omega)$, if $\omega$ is even.

**Main Results:**

**Lemma 1:** If G is $\{3K_1, K_{\omega+1}\}$-free with odd $\omega$, then $\chi \leq \frac{1}{8}(\omega^2+12\omega-13)$ where $\omega \leq 11$.

Proof: Let G be the smallest $\{3K_1, K_{\omega+1}\}$-free graph with $\chi(G) > \frac{1}{8}(\omega^2+12\omega-13)$ for odd $\omega \in \{3,..,11\}$. Then by minimality, $\forall$ u, $\chi(G-u) \leq \frac{1}{8}(\omega^2+12\omega-13)$ for $\omega \in \{3,..., 11\}$. Let deg u=$\Delta$. Then $\chi(G-u) = \frac{1}{8}(\omega^2+12\omega-13)$ and $\chi = \frac{1}{8}(\omega^2+12\omega-13) + 1$ for odd $\omega \in \{3,.., 11\}$.    **I**

Let C be a $\chi$-coloring of G in which u receives the unique color $\chi$. Let C = {1, 2, ....., $\chi$}. Clearly G has at the most two vertices with the same color. Let G have unique i-vertices for i ∈ {1,..,r} ⊆ C.

Let $V(G) = V_R \cup V_S \cup V_{S'} \cup V_T \cup V_{T'} \cup V_K$ where

$V_R$ = {v/ v receives color i, 1<=i<=r in C} with $|V_R|$= r. Clearly $<V_R>$ = $K_r$
$V_S$ = {v/ v is the only j-vertex in $<N(u)>$, j > r, j ∈ C}. Let $|V_S|$ = s.
$V_{S'}$ = {v/ v ∉ N(u) and v receives the same color as some w in $V_S$} = V(G) – N(u) and $|V_{S'}|$ = s.
$V_T$ = {v/ vw ∉ E(G) for some w ∈ $V_R$-u}. Let $|V_T|$ = t.
$V_{T'}$ = {v/ v receives the same color as some w in $V_T$}. Then $|V_{T'}|$ = t.

$V_K = V(G) - V_R - V_S - V_{S'} - V_T - V_{T'}$. Let $|V_K| = 2k$.

Clearly $r+s+t+k = \frac{1}{8}(\omega^2+12\omega-13)+1$ = Total no. of colors in C.  **II**

Also $\Delta = r-1+s+2t+2k$  **III**

Again as deg $u=\Delta$, and u is non-adjacent to s vertices in G, every vertex of G is non-adjacent to at least s vertices in G,  **IV**

**Claim 1:** Every $v \in V_R$ is adjacent to all vertices of $V_S \cup V_T$

If say $vw \notin E(G)$ for some $v \in V_R$ and $w \in V_S$, then color w by the color of v and u by color of w, a contradiction. Also if say $vw \notin E(G)$ for some $v \in V_R$ and $w \in V_T$, then $\exists\ w' \in V_{T'}$ and $z \in V_R$ s.t. $zw' \notin E(G)$. Then color w' by the color of z, w by the color of v and u by color of w, a contradiction.

This proves Claim 1.

**Claim 2:** Every $v \in V_S$ is adjacent to all vertices of $V_T$.
Let if possible $vw \notin E(G)$ for some $v \in V_S$ and $w \in V_T$. Let $w' \in V_{T'}$ have same color as w. Let $z \in V_R$ s.t. $zw' \notin E(G)$. Then color w' by the color of z, v by the color of w and u by color of v, a contradiction.

This proves Claim 2.

**Claim 3:** If $\omega<V_S>=m$ with $M \subseteq V_S$ containing the remaining s - m vertices of $V_S$, then every vertex of $V_R$-u is adjacent to at least $s-m+\omega<M>$ vertices of $V_{S'}$.

Let $s_1, s_2 \in V_S$ s.t. $s_1s_2 \notin E(G)$. Let $s_1', s_2' \in V_{S'}$ with same color as $s_1, s_2$ respectively. Let if possible $v \in V_R$ be non-adjacent to say $s_1'$, then color $s_1'$ by the color of v, $s_2$ by the color of $s_1$ and u by color of $s_2$, a contradiction. Thus if $s \in V_S$ is non-adjacent to some vertex of $V_S$, then every vertex of $V_R$ is adjacent to s. Now as deg $u = \Delta$, every vertex of G is non-adjacent to at least s vertices.

This proves Claim 3. This also implies that every vertex of $V_R$-u is non-adjacent to at least $t \geq s - m + \omega<M>$ vertices of $V_{T'}$.

**Claim 4:** If $\exists$ vertices $a', b' \in V_{T'}$ s.t. $a'u_i, b'u_j \notin E(G)$ for some $u_i, u_j \in V_R$- u, $i \neq j$ then $ab \in E(G)$ where a, b in $\in V_T$ have same colors as a', b' resply.

Let if possible $ab \notin E(G)$. Then color a' by i, b' by j, a by color of b and u by color of a, a contradiction.

Hence the claim 4 holds.

Thus $r+\omega<V_S>+\omega<V_T> \leq \omega$; $r + \omega<V_S \cup V_T \cup V_K> \leq \omega$  **V**

Next we prove the main results.

**1.    $\omega+1=4$**

As R(3, 3) = 6, $\Delta \leq 5$ and by [1], $\chi \leq 4 = \frac{1}{8}(\omega^2+12\omega-13)$

**2.    $\omega+1=6$**
Let if possible $\chi > 9$. As R(3, 6) = 18 and R(3, 5) = 14, $p \leq 17$ and $\Delta \leq 13$. Hence $p = 2\chi-r \Rightarrow 17 \geq$

20–r => r ≥ 3. Also r ≤ 5 as $<V_r> = K_r$. From **II**, r+s+t+k = 10 and from **III**, 13 ≥ Δ = r-1+s+2t+2k. Thus r+s ≥ 6.

As r ≤ 5, s ≥ 1, and s+t+k ≥ 5. If r ≥ 4, then by **V**, s+t+2k ≤ 2, a contradiction. Hence r=3 => s ≥ 3 and $\omega<V_S> \geq 2$. Then from **V**, t=0 and $\omega<V_S> \geq 2$. This implies that all vertices of $V_R$ are non-adjacent to all vertices of $V_S$, contrary to the Claim 3 as $\omega<V_S> < s$.

Hence $\chi \leq 9 = \frac{1}{8}(\omega^2+12\omega-13)$.

3.    ω+1=8

From [2] we get $\chi \leq \frac{1}{2}\lceil \Delta+\omega+1\rceil \leq \{(22+8)/2\} = 15 = \frac{1}{8}(\omega^2+12\omega-13)$. Hence $\chi \leq 15$.

Here in fact, we prove that $\chi \leq 14$. Let if possible $\chi = 15$. Then r+s+t+k = 15 and r, s ≤ 7. Also 27 ≥ p = 30–r => r ≥ 3 and 22 ≥ r-1+s+30-2r-2s => r+s ≥ 7. If r ≥ 4, then s+t+2k ≤ 8 and r ≥ 7. Thus r = 7 and s = t = k = 0. But then G ~ $K_7$ and $\chi = 7$, a contradiction. Hence r = 3, then s ≥ 4, s+t+k = 12 and s+t+2k ≤ 13 => k ≤ 1 and s+t ≥ 11. Thus 4 ≤ s ≤ 7 and t ≥ 4. Then $\omega<V_t> \geq 2$ and hence $\omega<V_s> = 2 = \omega<V_t>$. But then s, t ≤ 5, a contradiction.

Hence $\chi \leq 14$.

4.    ω+1=10
Let if possible $\chi > 22$. Here r+s+t+k = 23. By [3], R(3, 10) ≤ 42 and R(3, 9) = 36 => p ≤ 41 and Δ ≤ 35. 41 ≥ p = 2χ–r => r ≥ 5. Also r ≤ 9. Thus 5 ≤ r ≤ 9. Also 35 ≥ Δ = r-1+s+46–2r-2s => r+s ≥ 10. Again as r ≥ 5, by **V**, s+t+2k ≤ 13, a contradiction.

Hence $\chi \leq 22 = \frac{1}{8}(\omega^2+12\omega-13)$.

Here in fact, we prove that $\chi \leq 21$. Let if possible $\chi = 22$. Here r+s+t+k = 22, r, s ≤ 9. By [3], R(3, 10) ≤ 42 and R(3, 9) = 36 => p ≤ 41 and Δ ≤ 35. Thus 41 ≥ p = 2χ–r => r ≥ 3. Thus 3 ≤ r ≤ 9. Also 35 ≥ Δ = r-1+s+46–2r-2s => r+s ≥ 10. If r ≥ 4, then s+t+k ≤ s+t+2k ≤ 17 => r ≥ 5. But then again s+t+k ≤ s+t+2k ≤ 13 => r ≥ 9. Thus as before r = 9, s = t = k = 0, G ~ $K_9$ and χ = 9, a contradiction. Hence r = 3, s ≥ 7, s+t+k = 19 and s+t+2k ≤ 22. Thus k ≤ 3 and s+t ≥ 16. Now s ≤ 9 => t ≥ 7 and $\omega<V_t> \geq 3 => \omega<V_s> = 3 = \omega<V_t>$. Thus s, t ≤ 8 and in fact, s = t = 8. Clearly by Claim 4, every vertex of $V_r$ – u is non-adjacent to all vertices of $V_t$. Then by Claim 3, $\omega<V_t> = 8$, a contradiction.

Hence $\chi \leq 21$.

5.    ω+1 = 12
Let if possible $\chi > 30$. Then r+s+t+k = 31. By [3], R(3, 12) ≤ 59 and R(3, 11) ≤ 50 => p ≤ 58 and Δ ≤ 49. Thus 58 ≥ 62–r => r ≥ 4. Also r ≤ 11. Thus 20 ≤ s+t+k and by **V**, $\omega<V_S \cup V_T \cup V_K> \geq 6$ and r ≤ 5. Again as Δ = r-1+s+2t+2k, 49 ≥ r-1+s+62-2r-2s = 61-r–s => r+s >=12. Also s ≤ ω = 11.     **VI**

If r=5, then s+t+k = 26 and $\omega<V_S \cup V_T \cup V_K> \geq 7$, a contradiction. Hence r=4 and by **VI**, s >= 8 and $\omega<V_S> \geq 3$. Also s+t+k = 27 and by **V**, as $\omega<V_S \cup V_T \cup V_K> \leq 7$, k=0 and s+t=27. By **VI**, t ≥ 16 and $\omega<V_T> \geq 5$, a contradiction to **V**.

Hence $\chi \leq 30 = \frac{1}{8}(\omega^2+12\omega-13)$. This proves **Lemma 1**.

**Lemma 2:** If G is {$3K_1, K_{\omega+1}$}-free with even ω, then $\chi \leq \frac{1}{8}(\omega^2+10\omega)$, where ω ≤ 10

Proof: By [2], $\chi(G) \leq \lceil \frac{1}{2}(\Delta+\omega+1)/2\rceil \leq \lceil \frac{1}{2}(R(3, \omega)+\omega)\rceil$. For even ω with ω ≤ 8, it can be easily seen that $\lceil \frac{1}{2}(R(3, \omega)+\omega)\rceil = \frac{1}{8}(\omega^2+10\omega)$. Thus the result is true for even $\omega$ with $\omega \leq 8$.

Hence we prove the lemma for $\omega=10$ i.e. To prove that if G is $K_{11}$-free, then $\chi \leq 25$.

Let G be the smallest $\{3K_1, K_{11}\}$-free graph with $\chi(G) > 25$. Then $\chi(G-u) \leq 25$ $\forall$ u. Let deg u=$\Delta$. Thus $\chi(G-u)=25$ and $\chi=26$. **I**

Let C be a 26-coloring of G in which u receives the unique color 26. Let C = {1, 2, ....., 26}. Clearly G has at the most two vertices having the same color. Let G have unique i-vertices for i ∈ {1,..,r} ⊆ C.

Let $V(G) = V_R \cup V_S \cup V_{S'} \cup V_T \cup V_{T'} \cup V_K$ as defined in **Lemma 1** above. Then s ≤ 10.

Then r+s+t+k = 26 = Total no. of colors in C. **II**

Also by [3], 41 ≥ $\Delta$ = r-1+s+2t+2k = r-1+s+52–2r–2s => r+s ≥ 10 **III**

Again r+ $\omega<V_S>+\omega<V_T> \leq 10$; r+$\omega<V_S \cup V_T \cup V_K> \leq 10$ **IV**
Now 49 ≥ p = 52–r => r ≥ 3. Thus 3 ≤ r ≤ 10 => s+t+k ≥ 16 and $\omega<V_S \cup V_T \cup V_K> \geq 5$ => r ≤ 5.

Claims 1 to 4 of Lemma 1, hold good here also.

If r = 4, then s ≥ 6, s+t+k = 22 and $\omega<V_S+V_T+V_K> \leq 6$ => s+t+2k ≤ 22. Thus k=0, 6 ≤ s ≤ 10 and t ≥ 12. But then $\omega<V_S> \geq 3$ and $\omega<V_T> \geq 4$, a contradiction to **IV**.

If r = 3, then s ≥ 7, $\omega<V_S> \geq 3$, and s+t+k = 23 => $\omega<V_S+V_T+V_K> \geq 7$ and in fact $\omega<V_S+V_T+V_K>$ = 7. Thus s+t+2k ≤ 27, k ≤ 4, and s+t ≥ 19. As s ≤ 10, t ≥ 9 and $\omega<V_T> \geq 4$. Thus by **IV**, $\omega<V_S>$ = 3, s ≤ 8, $\omega<V_T>$ = 4 and 11 ≤ t ≤ 13. By Claim 3, every vertex of $V_r$-u is non-adjacent to at least 6 vertices of $V_{T'}$. Then as r=3, by Claim 4, it can be easily seen that $\omega<V_T> \geq 5$, a contradiction.

This proves Lemma 2.

Let G be a $\{3K_1, K_{\omega+1}\}$-free graph. Then the following table summarizes the above results:

| $\omega + 1 =$ | $\chi \leq$ |
|---|---|
| 3 | 3 |
| 4 | 4 |
| 5 | 7 |
| 6 | 9 |
| 7 | 12 |
| 8 | 14 |
| 9 | 18 |
| 10 | 21 |
| 11 | 25 |
| 12 | 30 |

**Table 1**

**Note:** The upper bounds for ω+1 = 5 and 6 were obtained in [13] using data from [12].

Next we make some conjectures.

**Conjecture 1:**
1. If ω is odd, then

$$\text{(A) } R(3, \omega) = \frac{1}{4}(\omega^2+8\omega-9), \text{ if } R(3, \omega) \text{ is even, and}$$

$$\text{(B) } R(3, \omega) = \frac{1}{4}(\omega^2+8\omega-13), \text{ if } R(3, \omega) \text{ is odd}$$

2. If ω is even, then

$$\text{(C) } R(3, \omega) = \frac{1}{4}(\omega^2+6\omega) \text{ if } R(3, \omega) \text{ is even and}$$

$$\text{(D) } R(3, \omega) = \frac{1}{4}(\omega^2+6\omega-4), \text{ if } R(3, \omega) \text{ is odd.}$$

The results can be verified for ω ≤ 9.

**Conjecture 2:** If G is $\{3K_1, K_{\omega+1}\}$-free and ω is odd, then $\chi \leq \lfloor \frac{1}{2}(\Delta+\omega+1) \rfloor$

The results can be easily verified for ω ≤ 9. If proved, then this is an improvement over Reed's conjecture for $\{3K_1, K_{\omega+1}\}$-free graphs with odd ω and Δ.

**Conjecture 3:** If G is $\{3K_1, K_{\omega+1}\}$-free,

$$3.1 \; \chi \leq \frac{1}{8}(\omega^2+12\omega-13), \text{ if } \omega \text{ is odd}$$

$$3.2 \; \chi \leq \frac{1}{8}(\omega^2+10\omega), \text{ if } \omega \text{ is even}$$